# DISCUSSION PAPER
# CONDITIONAL GROWTH CHARTS[1]

By Ying Wei and Xuming He

*Columbia University and University of Illinois*

Growth charts are often more informative when they are customized per subject, taking into account prior measurements and possibly other covariates of the subject. We study a global semiparametric quantile regression model that has the ability to estimate conditional quantiles without the usual distributional assumptions. The model can be estimated from longitudinal reference data with irregular measurement times and with some level of robustness against outliers, and it is also flexible for including covariate information. We propose a rank score test for large sample inference on covariates, and develop a new model assessment tool for longitudinal growth data. Our research indicates that the global model has the potential to be a very useful tool in conditional growth chart analysis.

**1. Introduction.** Growth charts, also known as reference centile charts, have been widely used to screen the measurements from an individual subject in the context of population values. A growth chart consists of a series of smooth curves plotted against time or another covariate, with each curve representing the trend of a given percentile of the measurement in a population. With several chosen percentile curves, a growth chart displays the distribution of a certain measurement within a certain range of time for a certain population. When a measurement is extreme on the chart, the subject is often identified for further investigation. An extreme measurement is likely to be a reflection of some unusual underlying physical condition.

Received December 2004; revised August 2005.
[1]Supported in part by National Science Foundation Grants DMS-01-02411 and DMS-05-04972, and National Security Agency Grant H98230-04-1-0034.
[2]Discussed in 10.1214/009053606000000632, 10.1214/009053606000000641, 10.1214/009053606000000650 and 10.1214/009053606000000669; rejoinder at 10.1214/009053606000000678.

The conventional method of constructing reference centiles is to get the empirical percentiles from cross-sectional data at a series of time points, and then fit a smooth polynomial curve to them. This method was used to develop the National Center of Health Statistics (NCHS) Growth Chart in 1977; see [15]. Cole [5] has become a classical work in the statistics literature on growth chart construction. In recent years, a number of different methods have been developed in the medical statistics literature; see, for example, [1, 9, 22, 32, 34, 39, 42]. Wright and Royston [43] have reviewed some of these methods. Most authors explore the beauty of normality to reduce quantile estimation to estimation of moments. Since most physical measurements are known to be nonnormal in distribution, a transformation to normality is generally used. The Box–Cox power transformation remains the most popular choice in this regard, as it can be found in a wide variety of medical journals—[7, 8, 29, 30], just to name a few. Arguably the most sophisticated transformation method is the LMS method of Cole [5]. It assumes that, at any time $t$, the measurement $Y(t)$ can be transformed to be approximately standard normal by $Z(t) = ((Y(t)/M(t))^{L(t)} - 1)/(L(t)S(t))$, where $L(t)$ is the best-fitting Box–Cox power for $Y(t)$, $M(t)$ is the median of $Y(t)$, and $S(t)$ is the scale of the transformed variable. The transformed variable $Z(t)$ is also known as the $z$-score. The penalized log-likelihood approach for estimating the $L$-$M$-$S$ functions has been proposed by Cole and Green [9].

Since conventional growth charts are usually developed from cross-sectional data, they are most useful for examining a subject with one measurement at a specific time. If a subject has more than one measurement, it can be more informative to study his growth path rather than a single measurement. Knowing a subject's prior growth path gives us a better understanding about his current growth status. For example, if a subject who has been growing along the 75th percentile curve suddenly drops below the median curve, it might be worth singling him out for further investigation. If we simply compare the growth path against the reference centiles generated from cross-sectional data, it would be hard to know how big a drop should be deemed unusual. Another important reason arguing against this approach is given in [6]. During both infancy and puberty, the subjects on the upper (or lower) centiles tend to grow toward the median at a faster growth rate than others, which is known as "catch-up" growth. Their growth paths are thus likely to move across the reference centiles during infancy or puberty. This is a normal phenomenon, but tracking those subjects on a conventional growth chart may give us an incorrect impression of abnormal growth. Clearly, one should compare a subject with the group of subjects with similar growth paths. In this paper, we focus on developing conditional growth charts based on longitudinal data.

Several authors, including Royston [33], have considered models for conditional reference centiles. Cole [6] considered the LMS-AR model $Z_t = a_t + b_t Z_{t-1} + e_t$, where $Z_t$ is the so-called $z$-score from a power transformation at time $t$, $Z_{t-1}$ is the lagged measurement at time $t-1$ and $e_t$ is distributed as $N(0, \sigma_t^2)$. In addition to the assumption of normality, this method requires fixed measurement time intervals. A more general model allowing varying measurement time intervals has been proposed by Thompson and Fatti [38], but they assumed a multivariate normal distribution for the measurements and the covariates at all time points and used the maximum likelihood estimator for the mean and variance functions. Scheike and Zhang [35] and Scheike, Zhang and Juul [36] considered a longitudinal regression model $Y_{i,j} = m(Y_{i,j-1}, t_{i,j} - t_{i,j-1}) + \varepsilon_{i,j}$ for height or weight, where $Y_{i,j}$ is the $j$th measurement of the $i$th subject at a random time $t_{i,j}$, $m(\cdot)$ is an unknown function of the prior measurement and the time duration from the prior measurement, and $\varepsilon_{i,j}$ is normally distributed. The reference centiles are constructed based on an estimate of $m$ and the normality of $Y_{i,j}$.

In practice, the longitudinal data collected for constructing reference centiles often exhibit the following features. First, the measurement time intervals are somewhat random. Even under the standard guideline of taking measurements on a nationally accepted schedule, say bimonthly for young children, the actual measurement times deviate from the fixed schedule due to practical considerations. A new subject to be screened for medical reasons is unlikely to follow the same schedule as the reference group. Second, transformation to normality is often a reasonable endeavor at a given time point, but the conditional models would require more than marginal normality of the $z$-scores. Normality of the error distributions in the conditional models can fail even if the marginal measurements are nearly normal; see an example later in the paper. Third, it is often desirable to incorporate subject-level covariates such as parent height in the conditional growth analysis. Normality in the conditional models becomes a stronger requirement in the presence of covariates. The existing approaches to conditional reference centiles can be ineffective in view of the data features discussed here. The semiparametric quantile regression model proposed in [40] aims to better accommodate data of this type; see illustrations in [41].

In this paper we further develop the global semiparametric approach to conditional growth charts by providing asymptotic theory, inference tools and a model assessment technique. Section 2 describes the global model and a simple estimation procedure. Section 3 provides the large sample properties of the estimates, and Section 4 gives a rank score test for covariate effects in the model. Section 5 compares the global model with the LMS-AR method of Cole [6] and illustrates the difficulty in preserving normality in conditional models. We also provide a new tool of model assessment for longitudinal growth data in Section 6. Some practical issues of using growth charts, both conditional and unconditional, are given in Section 7, and the technical proofs of our theorems are given in Section 8.

**2. A global model.** Suppose that we have $n$ subjects, and the $i$th subject has $m_i$ measurements at times $t_{i,1}, t_{i,2}, \ldots, t_{i,m_i}$, which are not necessarily evenly spaced. We denote as $Y_{i,j}$ the $j$th measurement of the $i$th subject, and $D_{i,j,k} = t_{i,j} - t_{i,j-k}$ is the time distance between the $j$th and $(j-k)$th measurements. For a given $\tau$, we model $Y_{i,j}$ by

$$(2.1) \qquad Y_{i,j} = g_0(t_{i,j}) + \sum_{k=1}^{p}(a_k + b_k D_{i,j,k})Y_{i,j-k} + X_{i,j}^{\top}\gamma_0 + e_{i,j}(\tau),$$

where $i = 1, \ldots, n$, $j = p+1, \ldots, m_i$, $X_{i,j} = (X_{i,j,1}, \ldots, X_{i,j,l})^{\top}$ consists of $l$ covariates for the $i$th subject at time $t_{i,j}$, and $e_{i,j}(\tau)$ is a random variable whose $\tau$th quantile, given the growth path up to the $(j-1)$st measurement and $X_{i,j}$, is zero. The conditional quantile of $Y_{i,j}$ given the $p$ prior measurements and the covariate $X_{i,j}$ includes $g_0$ as a nonparametric intercept function of the current measurement time, an autoregressive function of $Y_{i,j-1}, \ldots, Y_{i,j-p}$ whose coefficients are linear functions of measurement time distances $D_{i,j,k}$, and a linear function of the covariate $X_{i,j}$.

In general, the subjects take measurements at irregular time intervals. It is natural to assume that the dependence between two measurements varies with their measurement time distance. This motivates the choice of the autoregressive coefficients as functions of $D_{i,j,k}$. If all the subjects take measurements at the fixed time intervals, the coefficients simplify to constants. On the other hand, it is possible to generalize the model by using a more general form as considered by Cai and Xu [2]. We defer this discussion to Section 7.

We have intentionally written the error term in model (2.1) as $e_{i,j}(\tau)$, which is $\tau$-specific. This helps distinguish the quantile regression model from most other models of the same form. Model (2.1) allows the error term to be dependent on the covariates as well, so $e_{i,j}(\tau)$ is merely a convenient notation for the difference between $Y_{i,j}$ and its $\tau$th conditional quantile function. If the error terms are independent of $\tau$ and of the covariates in the model, then models of this form have been well studied in the longitudinal data analysis by many other authors, including Chiang, Rice and Wu [4], Lin and Carroll [24], Lin and Ying [23] and Fan and Li [11], to mention just a few.

For the sake of simpler presentation in the rest of the paper, we will convert to the usual notation of writing $e_{i,j}(\tau)$ as $e_{i,j}$, but it is always helpful to keep in mind that this error term is $\tau$-specific.

In (2.1), $g_0(t)$ is a nonparametric function with a certain degree of smoothness. Suppose that the range of interest in time is $t \in [t_L, t_U]$. Different smoothing methods may be used to estimate $g_0$, but Wei [40] used the convenient approach of regression splines. Specifically, we approximate $g_0(t)$ by a linear combination of B-spline basis functions as in [18]. Given a set of knots, we denote as $\pi(t) = (\pi_1(t), \pi_2(t), \ldots, \pi_{k_n}(t))^\top$ the set of $k_n$ B-spline basis functions of, say, order 4 (cubic). Let $\pi_{i,j} = \pi(t_{i,j})$, $\tilde{X}_{i,j} = (Y_{i,j-1}, D_{i,j,1} Y_{i,j-1}, \ldots, Y_{i,j-p}, D_{i,j,p} Y_{i,j-p}, X_{i,j}^\top)^\top$; we can then use $\pi_{i,j}^\top \hat{\boldsymbol{\alpha}}_n + \tilde{X}_{i,j}^\top \hat{\boldsymbol{\beta}}_n$ to estimate the $\tau$th conditional quantile function of $Y_{i,j}$, where $\hat{\boldsymbol{\alpha}}_n$ and $\hat{\boldsymbol{\beta}}_n$ are obtained by minimizing

$$(2.2) \qquad \sum_{i=1}^{n} \sum_{j=p+1}^{m_i} \rho_\tau(Y_{i,j} - \pi_{i,j}^\top \boldsymbol{\alpha} - \tilde{X}_{i,j}^\top \boldsymbol{\beta}),$$

where $\rho_\tau(r) = r(\tau - I(r < 0))$ is the quantile objective function; we refer to [26] for details.

Obviously, other smoothing methods (e.g., kernel smoothing) may be used in lieu of regression splines. One advantage of spline smoothing here is that in typical applications it is usually sufficient to preselect a set of knots in the interval $(t_L, t_U)$ using our general understanding of growth patterns. For example, it would be useful to place more knots during infancy and puberty than during other times. In this paper, we do not go into the issue of automated knot selection.

Optimization of (2.2) can be performed efficiently by linear programming techniques. In fact, the solution to (2.2) only requires software for quantile regression in linear models. The R package *quantreg* contributed by Koenker and the SAS PROC QUANTREG are handy for obtaining the estimates of $\boldsymbol{\alpha}$ and $\boldsymbol{\beta}$.

We further note that in model (2.1), even if the $e_{i,j}$'s are assumed to be independent, the within-subject correlation may be captured by the autoregressive component of the model. For a consistent estimate of the conditional quantile function, we do not need the i.i.d. assumption on $e_{i,j}$. The inference tools developed later in this paper are proven under stronger assumptions but with good robustness against common deviations from the i.i.d. assumption.

**3. Large sample property.** Under appropriate conditions, the intercept function estimate and the coefficient estimate for $\boldsymbol{\beta}$ converge to their true values at the optimal rates. To state the conditions, we use $\|\cdot\|$ for the Euclidean norm, and $a_n \approx b_n$ to mean $0 < \liminf_n a_n/b_n < \limsup_n a_n/b_n < \infty$. Throughout, we assume that the number of measurements per subject is bounded, but the number of subjects grows. For the sake of convenience and without loss of generality, we assume that $p = 1$ in model (2.1), and consequently drop the subscript $k$ from $D_{i,j,k}$ so that $D_{i,j} = D_{i,j,1}$. Also let $\psi(x) = \rho'_\tau(x) = \tau - I_{\{x \leq 0\}}$, and $\psi(e_i) = (\psi(e_{i,2}), \psi(e_{i,3}), \ldots, \psi(e_{i,m_i}))^\top$. Since our results will be stated for a given $\tau$, we have dropped the dependence of $\psi$ on $\tau$ here.

Following He and Shi [18], we assume a spline approximation to $g_0(t)$ as $\pi_{i,j}^\top \boldsymbol{\alpha}_0$, with $R_{nij} = \pi_{i,j}^\top \boldsymbol{\alpha}_0 - g_0(t_{i,j})$ as the approximation error. We state the assumptions for this section.

ASSUMPTION 1. $g_0(t)$ has bounded $r$th derivative for some $r \geq 1$.

ASSUMPTION 2. The numbers of measurements $m_i$ are uniformly bounded for all $n$.

ASSUMPTION 3. The errors $e_i = (e_{i,2}, \ldots, e_{i,m_i})^\top$ have unspecified variance–covariance structures but are independent across subjects.

ASSUMPTION 4. $A_i = E[\psi(e_i)\psi(e_i)^\top] > 0$ for each $i$ and $\sup_i \|A_i\| < \infty$.

Let $b_{i,j} = f_{e_{i,j}}(0)$, where $f_{e_{i,j}}$ is the conditional p.d.f. of $e_{i,j}$ given $t_{i,j}$ and $\tilde{X}_{i,j}$.

ASSUMPTION 5. $0 < \inf_{i,j} b_{i,j} < \sup_{i,j} b_{i,j} < \infty$, and as $s \to 0$, $\sup_{i,j} |E\psi(e_{i,j} + s) - b_{i,j}s| = O(s^2)$.

ASSUMPTION 6. There exists an $(l+2)$-dimensional function $\eta(t) = (\eta_1(t), \eta_2(t), \ldots, \eta_{l+2}(t))$ with bounded $r$th derivative such that for any $i$ and $j$,

$$\tilde{X}_{i,j} = \eta(t_{i,j}) + \delta_{i,j},$$

where $\delta_{i,j} = (\delta_{i,j,1}, \delta_{i,j,2}, \ldots, \delta_{i,j,l+2})$ are random vectors with mean zero and are independent of $e_{i,j}$ given $(i,j)$. Let $\tilde{X}_n = (\tilde{X}_{1,2}, \ldots, \tilde{X}_{1,m_1}, \tilde{X}_{2,2}, \ldots, \tilde{X}_{n,m_n})$ be the $(l+2) \times N$ matrix with $N = \sum_{i=1}^n (m_i - 1)$, and let $\Delta_n$ and $V_n$ be defined in a similar way using $\delta_{i,j}$ and $\eta(t_{i,j})$, respectively. Then we have $\tilde{X}_n = V_n + \Delta_n$, where $V_n$ is the mean of $\tilde{X}_n$ so that $E\Delta_n = 0$.

ASSUMPTION 7. $\sup_n n^{-1} E\|\Delta_n\|^2 < \infty$.

Assumption 7 is readily satisfied if all the variables have finite supports. Note that for any $i$, $(t_{i,1}, t_{i,2}, \ldots, t_{i,m_i})$ are order statistics of the time variable, thus $t_{i,j-1}$ and $Y_{i,j-1}$ are not independent of $t_{i,j}$. Therefore, $\eta_1(t_{i,j})$ is the conditional mean of the prior measurement given the current time, and $\eta_2(t_{i,j})$ is the conditional mean of $D_{i,j}$ times $Y_{i,j-1}$ given current time $t_{i,j}$. Similarly, $(\eta_3, \ldots, \eta_{l+2})$ are the conditional means of the covariates $X_{i,j}$ given $t_{i,j}$.

Let $Z_n = (\pi_{1,2}, \pi_{1,3}, \ldots, \pi_{i,j}, \ldots, \pi_{n,m_n})_{N \times k_n}^\top$, $G_n = Z_n(Z_n^\top Z_n)Z_n^\top$ and $B_n = \mathrm{diag}(b_{i,j})$, where $\{b_{i,j}\}$ is sequenced by indices $i$ and $j$ as in $Z_n$. For additional assumptions, we define $X_n^* = (I - G_n)\tilde{X}_n$ and $K_n = X_n^{*T} B_n X_n^*$. If we denote as $X_n^{*(i)}$ the $m_i \times (l+2)$ submatrix of $X_n^*$ corresponding to the $i$th subject, and $B^{(i)} = \mathrm{diag}(b_{i,2}, \ldots, b_{i,m_i})$, then $K_n = \sum_{i=1}^n X_n^{*(i)T} B^{(i)} X_n^{*(i)}$. By the assumption of between-subject independence, $K_n$ can be regarded as an independent sum. Furthermore, let $\Lambda_n = \mathrm{diag}(A_1, \ldots, A_n)$ and $S_n = X_n^{*T} \Lambda_n X_n^*$. Then we have

ASSUMPTION 8. *There exists a positive definite matrix $K$ such that $K_n/n \to K$ in probability, and $E(K_n/n) \to K$.*

ASSUMPTION 9. *There exists a nonnegative definite matrix $S$ such that $S_n/n \to S$ in probability.*

Let $\hat{g}_n(t_{i,j}) = \pi_{i,j}^\top \hat{\boldsymbol{\alpha}}_n$; the following theorem summarizes the large sample properties of $\hat{\boldsymbol{\beta}}_n$ and $\hat{g}_n$.

THEOREM 3.1. (i) *Under Assumptions 1–7, if the number of knots $k_n \approx n^{1/(2r+1)}$, then*

$$\|\hat{\boldsymbol{\beta}}_n - \boldsymbol{\beta}_0\| = O_p(n^{-r/(2r+1)}) \tag{3.1}$$

*and*

$$\frac{1}{n} \sum_{i=1}^n \sum_{j=p+1}^{m_i} (\hat{g}_n(t_{i,j}) - g_0(t_{i,j}))^2 = O_p(n^{-2r/(2r+1)}). \tag{3.2}$$

(ii) *With Assumptions 8 and 9 additionally, we have*

$$\sqrt{n}(\hat{\boldsymbol{\beta}}_n - \boldsymbol{\beta}_0) \to N(0, \Sigma), \tag{3.3}$$

*where $\Sigma = K^{-1} S K^{-1}$.*

**4. Rank score test.** The asymptotic normality of $\hat{\boldsymbol{\beta}}_n$ enables us to make large sample inference on the coefficients $a_k$, $b_k$ and $\gamma_0$. However, the variance–covariance matrix $\Sigma$ of Theorem 3.1 is not easy to estimate, as it involves the densities of $e_{i,j}$. Likelihood-based tests such as those discussed in [10] are harder to develop for the quantile models than for conditional mean regression. The rank score test given by Gutenbrunner, Jurečková, Koenker and Portnoy [13] provided an attractive alternative to hypothesis testing on quantile regression coefficients by avoiding direct estimation of the error densities. In this section we extend the rank score test to the semiparametric model (2.1). Fundamental theories of rank-based inference can be found in [12, 14, 28].

As in [13], we assume $e_{i,j}$ to be i.i.d. in this section. Even so, the inclusion of a nonparametric component $g_0$ in the model complicates the derivation of the limiting distribution of a rank score test statistic. Fortunately, the $\chi^2$ limiting distribution initially derived from linear quantile regression models remains valid when $g_0$ is estimated by regression splines.

4.1. *Test statistic and asymptotic distribution.* Suppose that $\boldsymbol{\beta}_{01}$ is a $q$-dimensional subset of $\boldsymbol{\beta}_0$ in (2.1), and $\tilde{X}_{n1}^\top$ is the $N \times q$ design matrix corresponding to $\boldsymbol{\beta}_{01}$. Let $\boldsymbol{\beta}_{02}$ and $\tilde{X}_{n2}$ be the remaining parameter vector and design matrix, respectively; then the global model (2.1) has a partitioned form

$$(4.1) \qquad Y_n = g_0(T) + \tilde{X}_{n1}^\top \boldsymbol{\beta}_{01} + \tilde{X}_{n2}^\top \boldsymbol{\beta}_{02} + e_n.$$

Suppose that we wish to test the null hypothesis $H_0 : \boldsymbol{\beta}_{01} = 0$, versus the alternative hypothesis $H_1 : \boldsymbol{\beta}_{01} \neq 0$. Using the spline approximation described in Section 2, we further denote that

$$\begin{aligned}
Z_n &= (\pi_{1,2}^\top, \pi_{1,3}^\top, \ldots, \pi_{1,m_1}^\top, \ldots, \pi_{n,2}^\top, \pi_{n,3}^\top, \ldots, \pi_{n,m_n}^\top)^\top, \\
\mathcal{W}_n &= (Z_n^\top, \tilde{X}_{n2}^\top)^\top, \\
G_n &= \mathcal{W}_n^\top (\mathcal{W}_n^\top \mathcal{W}_n)^{-1} \mathcal{W}_n, \\
\mathcal{V}_n^\top &= G_n \tilde{X}_{n1}^\top, \\
\phi_0 &= (\alpha_0^\top, \boldsymbol{\beta}_{02}^\top)^\top,
\end{aligned}$$

where $\alpha_0$ is defined in Section 3 and $\mathcal{W}_n$ is the pseudodesign matrix when $H_0$ is true.

Let $v_{i,j}$ be the column component of $\mathcal{V}_n$ that corresponds to the $i$th subject and the $j$th measurement. Similarly, we denote $w_{i,j}$ as the column component of $\mathcal{W}_n$ that corresponds to the $i$th subject and the $j$th measurement. The dimension of $w_{i,j}$ is increasing in the dimension of the B-spline space.

Now, let $S_n = n^{-1/2} \sum_{i=1}^n \sum_{j=2}^{m_i} \varphi_\tau(Y_{i,j} - w_{i,j}^\top \hat{\phi}_n) v_{i,j}$, where $\hat{\phi}_n = \arg\min_\phi \sum_{i=1}^n \sum_{j=2}^{m_i} \rho_\tau(Y_{i,j} - w_{i,j}^\top \phi)$ and $V_n = n^{-1} \tau(1-\tau) \tilde{X}_{n1} (I_n - G_n) \tilde{X}_{n1}^\top$. We define the rank score test statistics as

$$(4.2) \qquad T_n = S_n^\top V_n^{-1} S_n.$$

To obtain the asymptotic distribution of $T_n$, we make additional assumptions as follows.

ASSUMPTION D1. The random errors $e_{i,j}$ are independent of one another, and there exists a constant $b$ such that $b_{i,j} = f_{e_{i,j}}(0) = b$, for all $i$ and $j$.

ASSUMPTION D2. $\max_{i,j} \|w_{i,j}\|^2 = O_p(k_n)$, $E[\max_{i,j} \|v_{i,j}\|^2] < \infty$ and $\sup_{\|\alpha\| \leq 1} \sum_{i=1}^{n} \sum_{j=2}^{m_i} \|v_{i,j}\|^2 \|w_{i,j}^\top \alpha\|^2 = O_p(n)$.

ASSUMPTION D3. There exists a constant $\kappa$ such that $\sup_{i,j} |f_{Y_{i,j}|v_{i,j},w_{i,j}}| \leq \kappa$.

In Assumption D3, we assume the conditional density of $Y_{i,j}$ to be uniformly bounded. Assumption D1 holds if the $e_{i,j}$'s are i.i.d. Besides Assumptions D1–D3, we have a slightly stronger assumption on the smoothness of the function $g_0(t)$. We modify Assumption 1 in Section 3 as

ASSUMPTION 1*. $g_0(t)$ has bounded $r$th derivative for some $r > 2$.

We now have the main theorem of this section, which allows, but is not limited to, the choice of $k_n \approx n^{1/(2r+1)}$ given in Theorem 3.1. Using $a_n \ll b_n$ to mean $\lim_{n \to \infty} a_n/b_n = 0$, we have the following result.

THEOREM 4.1. *Under Assumptions 1*, 2–7 and D1–D3, if the number of knots $k_n$ for the B-spline space satisfies $n^{1/4r} \ll k_n \ll n^{1/4}$ with $r > 2$, we have $T_n \to \chi_q^2$ as $n \to \infty$.*

It is helpful to note that the limiting distribution of $T_n$ is invariant over a wide range of choices for $k_n$. Practically speaking, the result suggests that the number of knots is not a critical issue for the rank score tests. In our empirical investigations, the validity of the test appeared robust against violation of the i.i.d. assumption on $e_{i,j}$.

4.2. *An assessment with Monte Carlo simulations.* We may assess the performance of the rank score test described in the preceding subsection by Monte Carlo. A set of models is selected with different error distributions, percentile levels $\tau$ and model structures. In each model, the number of subjects is 200, the number of measurements from each subject is 10 and the percentile level is 0.5 or 0.95. The number of Monte Carlo samples used for estimating the Type I error is 10,000 in each scenario, so that the standard error of the estimate is around 0.2%.

MODEL 1. *A semiparametric longitudinal model with i.i.d. random errors.* $y_{i,j} = h(t_{i,j}) + by_{i,j-1} + x_i + e_{i,j}$, where $y_{i,j}$ is the $j$th measurement from the $i$th subject at time $t_{i,j}$, $h(t) = 40t/(1 + 4t)$ is the intercept function, $e_{i,j}$ are i.i.d. standard normal, $x_i$, independent of $e_{i,j}$, are normal with mean 10 and standard deviation 1, and $b$ is a parameter of interest. Here, $h(t)$ has similar shape and scale as a reference centile curve of infant weight. Without loss of generality, we assume that, given $i$, the measurement times $(t_{i,1}, t_{i,2}, \ldots, t_{i,10})$ are order statistics of ten uniformly distributed random variables on the interval $[0,1]$. The intercept function $h(t)$ is monotone increasing and is approximated by a cubic spline with internal knots (0.25, 0.5, 0.75). Using the rank score test, we test the null hypothesis: $H_0 : b = 0$, that is, the significance of the autoregressor is of interest. Models 2 and 3 have the same structure as Model 1, but assume that the distribution of $e_{i,j}$ is $t_3/\sigma_{t_3}$ and $\chi_1^2/\sigma_{\chi_1^2}$, respectively, where the variances of $e_{i,j}$'s are standardized to 1.

TABLE 1
*Type* I *errors for Models* 1–6

|  | $\tau = 0.5$ | | | $\tau = 0.95$ | | |
| --- | --- | --- | --- | --- | --- | --- |
| # of knots | Model 1 | Model 2 | Model 3 | Model 1 | Model 2 | Model 3 |
| 1 | 0.0524 | 0.0535 | 0.0531 | 0.0508 | 0.0572 | 0.0523 |
| 2 | 0.0494 | 0.0472 | 0.0488 | 0.0498 | 0.0551 | 0.0505 |
| 3 | 0.0513 | 0.0485 | 0.0479 | 0.0495 | 0.0550 | 0.0500 |
|  | **Model 4** | **Model 5** | **Model 6** | **Model 4** | **Model 5** | **Model 6** |
| 1 | 0.0534 | 0.0497 | 0.0563 | 0.0543 | 0.0672 | 0.0521 |
| 2 | 0.0506 | 0.0462 | 0.0505 | 0.0504 | 0.0662 | 0.0517 |
| 3 | 0.0510 | 0.0457 | 0.0523 | 0.0479 | 0.0684 | 0.0510 |

MODELS 4–6. *A longitudinal model with non-i.i.d. errors.*

$$y_{i,j} = h(t_{i,j}) + by_{i,j-1} + x_i + (1/2 + y_{i,j-1}/25)e_{i,j}.$$

For these models, the conditional quantile function is a linear function of $y_{i,j-1}$ with the autoregressive coefficient as $b_\tau = b + Q_\tau(e_{i,j})/25$, where $Q_\tau(e_{i,j})$ is the $\tau$th quantile of $e_{i,j}$. With Models 4–6, we test the null hypothesis $H_0 : b_\tau = 0$.

In all the models, the intercept function $h(t)$ is approximated by a regression cubic spline, which may vary with different knot selections. To see whether the rank score test is sensitive to the number of knots, we estimate $h(t)$ in Models 1–6 with one or two or three uniformly spaced knots. Table 1 gives the estimated Type I errors when $\tau = 0.5$ and $\tau = 0.95$ for the two sets of models. All of them are close to 0.05, the desired significance level.

A look at the estimated power curves (not presented in this paper) confirms that the performance of the test does not depend heavily on the choice of $k_n$ and is quite robust against modest deviations from the i.i.d. assumptions on the errors.

**5. Example on infant weight.** In this section we use the Finnish growth data [31] to demonstrate the use of the global model for subject screening and compare the results with the LMS-AR approach of Cole [6]. The data we use here include measurements from a total of 1143 Finnish boys from birth to 2 years of age.

5.1. *Global model on infant weight.* Monitoring the growth of weight during infancy is clinically important. For example, there is evidence that the rapid rates of weight gain during infancy could lead to obesity later in the childhood and may also be related to cardiovascular diseases later in life [37]. We use weight as the measurement for screening, but two prior measurements on weight and the current height are included in the model. Approximately, the measurements were taken monthly during this period of infancy, but few children were measured exactly on a monthly schedule.

Let $W_{i,j}$ and $H_{i,j}$ be the $j$th weight and height, respectively, of the $i$th subject at age $t_{i,j}$. The global model (2.1) for the $\tau$th quantile can be further specified as

$$(5.1) \quad W_{i,j} = g_\tau(t_{i,j}) + \sum_{k=1}^{2}(a_{k,\tau} + b_{k,\tau}D_{i,j,k})W_{i,j-k} + c_\tau H_{i,j} + e_{i,j}.$$

Note that the error variables $e_{i,j}$ are $\tau$-dependent in the model. A set of quantiles with $\tau \in \{0.03, 0.1, 0.25, 0.5, 0.75, 0.9, 0.97\}$ is chosen for consideration. In the example, we use cubic splines with evenly spaced internal knots (0.5, 1, 1.5).

Table 2 lists the estimates of $a_{k,\tau}$, $b_{k,\tau}$ ($k=1,2$) and $c_\tau$ for the seven $\tau$'s. The numbers in parentheses are the corresponding $p$-values based on the rank score test. They suggest that the second prior weight is not significant at several $\tau$ levels. Even when it is statistically significant, the small magnitudes of the coefficients mean that the contribution to the conditional quantiles would be small relative to that from the first prior weight. Height is a significant factor in the model.

Also note that $a_{1,\tau}$ decreases, but $b_{1,\tau}$ and $c_\tau$ increase with $\tau$. This means that the measurement time distance and height play larger roles for heavier boys than for lighter boys.

A boy's growth status can be assessed by comparing his current weight with a set of conditional quantiles estimated from model (5.1). By bootstrapping subjects, we can also get confidence intervals on the conditional quantiles.

The left panel in Figure 1 displays the weight growth path of one subject in the Finnish growth data. The growth path is plotted against the conventional (unconditional) centile curves of weight. We note that the weight of this subject had been growing along the 0.25th quantile but jumped to about the median level at age 0.61 year. A rapid weight gain could be a warning signal for obesity. Therefore, it is worthwhile to single him out early at age 0.61.

TABLE 2
*Estimated parameters for infant weight*

|  | $\tau$ | | | | | | |
|---|---|---|---|---|---|---|---|
|  | 0.03 | 0.10 | 0.25 | 0.50 | 0.75 | 0.90 | 0.97 |
| $a_{1,\tau}$ | 0.787 | 0.823 | 0.810 | 0.801 | 0.737 | 0.612 | 0.482 |
|  | (0.000) | (0.000) | (0.000) | (0.000) | (0.000) | (0.000) | (0.000) |
| $b_{1,\tau}$ | 0.163 | 0.216 | 0.279 | 0.308 | 0.332 | 0.334 | 0.346 |
|  | (0.000) | (0.000) | (0.000) | (0.000) | (0.000) | (0.000) | (0.000) |
| $a_{2,\tau}$ | 0.012 | −0.015 | −0.016 | −0.032 | −0.048 | −0.056 | −0.037 |
|  | (0.496) | (0.103) | (0.015) | (0.000) | (0.000) | (0.002) | (0.295) |
| $b_{2,\tau}$ | 0.045 | 0.020 | 0.008 | 0.023 | 0.021 | 0.043 | 0.071 |
|  | (0.004) | (0.083) | (0.547) | (0.080) | (0.223) | (0.012) | (0.040) |
| $c_\tau$ | 0.038 | 0.043 | 0.051 | 0.059 | 0.080 | 0.109 | 0.138 |
|  | (0.000) | (0.000) | (0.000) | (0.000) | (0.000) | (0.000) | (0.000) |

The numbers in the parentheses are $p$-values.

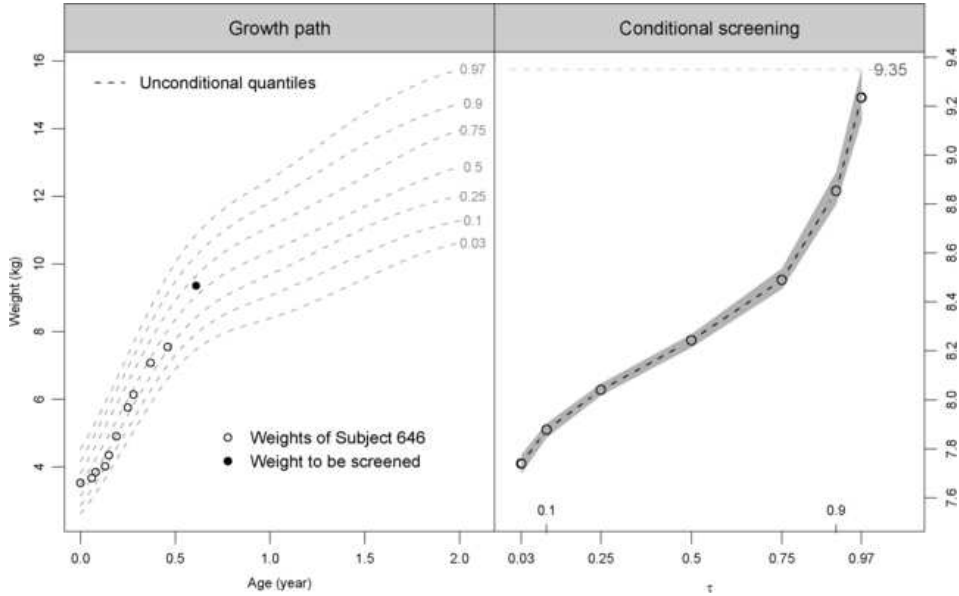

FIG. 1. *An example of conditional screening based on the global model.*

We apply model (5.1) to screen the weight at 0.61 year as the "current" measurement age. The right panel in Figure 1 shows the screening results from the global model, with the circle dots as the estimated conditional quantiles at age 0.61 from the global model. The boy's current weight is 9.35 kg, which is higher than the 0.97th conditional quantile, suggesting that he might be overweight given his prior path and current height. The gray band represents the pointwise 90% confidence band for the conditional quantiles obtained from the bootstrap method. The subject's current height is higher than the upper bound of the 0.97th conditional quantile, which reinforced our conclusion that, from ages 0.46 to 0.61, the boy had been gaining weight at an excessive rate, even though his weight was still below the 0.75th quantile at age 0.61 based on the conventional growth chart.

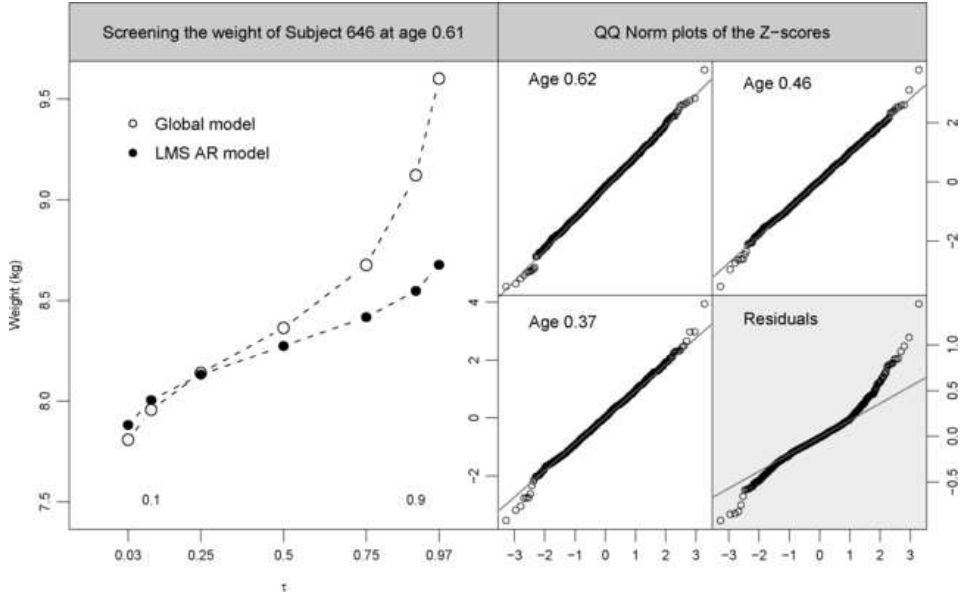

Fig. 2. *A comparison between the global model and the LMS-AR model.*

5.2. *Comparison with the LMS-AR model.* Cole [6] suggested using the LMS method to transform data to normality and then applying an autoregressive model to the $z$-scores to estimate the conditional quantiles. To use the LMS-AR model, we interpolate the data to have measurements on fixed time intervals. More specifically, we first transform the infant weight using the LMS method with the effective dimensions 7, 10 and 7 for the functions $L$, $M$ and $S$, respectively. We denote as $z_t$ the $z$-score at age $t$. An AR(2) model is then used at 0.61 year, with the first two lags at 0.46 year and 0.37 year as autoregressors:

$$(5.2) \qquad z_{0.61} = \alpha_0 + \alpha_1 z_{0.46} + \alpha_2 z_{0.37} + e,$$

where $e \sim N(0, \sigma^2)$.

The estimated quantiles from the LMS-AR(2) model, together with those from the global model using the same covariates, are plotted in the left panel of Figure 2. The two sets of conditional quantiles agree fairly well at the lower percentiles, although the former are slightly higher than the estimates from the global model. However, for the upper quantiles we observe considerably large differences between the two estimates: the estimates from the LMS-AR(2) model are much lower than the ones from the global model, and the difference widens as $\tau$ increases.

To understand the differences, we note that the LMS-AR approach implicitly assumes joint normality of the $z$-scores at three measurement times. The right panel of Figure 2 gives the normal QQ-plots of the $z$-scores for these four variables. While the LMS transformation does a pretty good job in achieving marginal normality at these three measurement times, it is clear that the distribution of the residuals deviates substantially from normality, and thus the joint normality is not achieved. The QQ-plot of the residuals shows clearly a heavy upper tail, which explains why the LMS-AR(2) method underestimates the upper quantiles.

**6. Model diagnosis.** To use the global model for subject screening, one has to have some faith in the adequacy of the global model. Does the model fit the data well? Can a more parsimonious model fit the data almost as well? These questions lead us to consider model assessment tools for the conditional model.

Goodness-of-fit tests developed for linear quantile regression models by Koenker and Xiao [27] and He and Zhu [20] may be extended to semiparametric models, but they are not made for model assessment over a set of quantiles. In this section, we propose an assessment tool tailored to our conditional models by trying to compare the empirical distribution of $Y$ with the simulated distribution from the model.

Because of the irregularity in measurement times, we would not have enough data for model assessment at any given time. However, we view $Y_{i,j}$, the $j$th measurement of the $i$th subject, as a random variable of interest for any fixed $j$. We allow $t_{i,j}$ to be random, and consider the distribution of the $j$th measurement for any group of subjects.

Given the quantile models, there is a nice way to generate $Y_{i,j}$ from its conditional distribution. Specifically, we draw $U$ from the uniform distribution on $(0, 1)$, and then generate $Y$ as the $U$th quantile from the (estimated) model. A resulting random sample follows the model-based conditional distribution.

For a given $j$, we randomly draw a subject from those who have at least $j$ measurements. Given this subject's $Y(t_{j-1})$ and $X(t_{j-1})$, we draw $Y(t_j)$ from the quantile models as described above. When this process is repeated many times, we obtain a simulated sample for the $j$th measurement. The size of the simulated sample could be as large as we wish, but 5000 is used in our example below. If the model fits the data well, the distribution of the simulated $Y(t_j)$ should match the distribution of the observed $j$th measurements. It is a weak requirement that the two marginal distributions should match; instead we would expect a good match for any subgroup of subjects defined by the covariates in the model, which forms the basis of our model assessment tool.

To have a quick idea of how well the empirical and the simulated distributions match, we may use the familiar QQ-plot. Another useful option is to look at $\sqrt{n}(\hat{\tau} - \tau)/(\tau(1-\tau))^{1/2}$ for a set of $\tau$'s, where $\hat{\tau}$ is the proportion of the observed $Y_{i,j}$ that fall below the $\tau$th quantile of the simulated distribution, and $n$ is the number of the observed $Y_{i,j}$.

For the infant weight example in the previous section, we display the QQ-plots in Figure 3 for $j = 2, 4, 6$ and 8. The reference line in each plot is $y = x$. These QQ-plots suggest a good match between the two distributions. The upper half of the figure plots the standardized $\hat{\tau} - \tau$ for 10 equally spaced $\tau$'s from 0.05 to 0.95. No differences are worrisomely large, except that the match at $j = 2$ was not as satisfactory as the others.

However, it is more meaningful to assess the agreement for subgroups of subjects. When all the subjects are included, we are only concerned with the marginal distributions, so even an unconditional model would do well.

We now consider some subgroups based on birth weight or current height and compare three models: the unconditional model, the global model (5.1) and a smaller global model obtained by dropping the second lagged term as well as height. Figure 4 gives the QQ-plots and the standardized $(\hat{\tau} - \tau)$ plots for two subgroups: Group 1 subjects whose birth weights are below the first quartile, and Group 2 subjects whose current heights are below the first quartile. Here, model adequacy is evidently different. The unconditional model, as expected, would overestimate the lower quantiles for Group 1 and Group 2. The smaller conditional model would fit the data better than the unconditional model, but the global model (5.1) does much better. The improvement due to lag-2 weight and due to height information is clearly demonstrated through the model assessment plots proposed in this section.

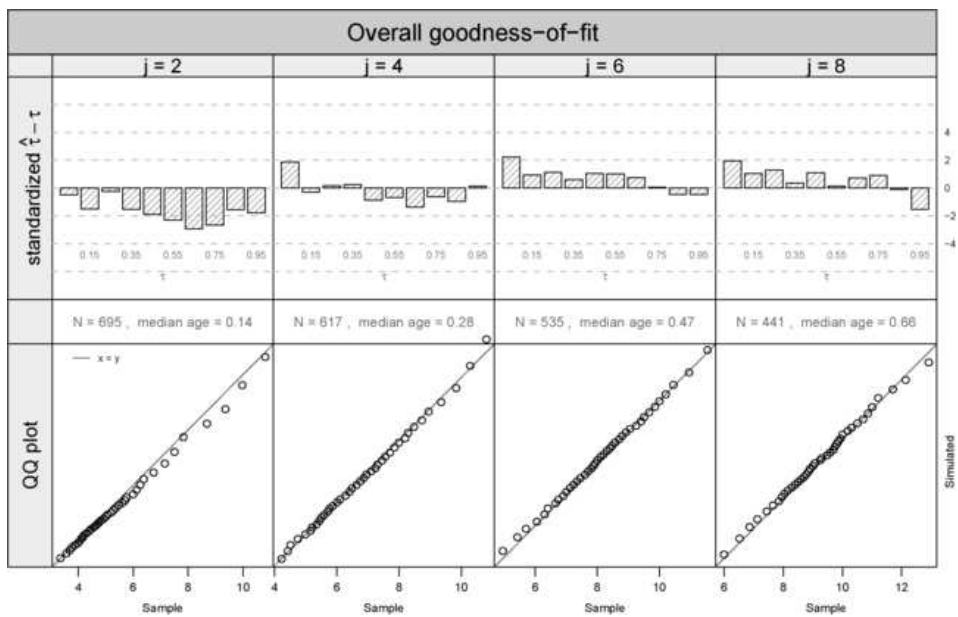

FIG. 3. *Checking goodness-of-fit at the $j$th measurement, where $N$ stands for the number of observations whose median measurement age (in year) is also displayed.*

By examining the model checking plots for various subgroups defined by covariates, we find that global model (5.1) is a good fit to the Finnish infant data for the purpose of conditional weight screening. Here, we choose not to specify a default for formations of subgroups; rather we leave it to the individual researcher to slice and detect. The model assessment ideas presented here are applicable to conditional quantile models in other applications, and formal lack-of-fit tests may be developed from these ideas. We defer them to future work.

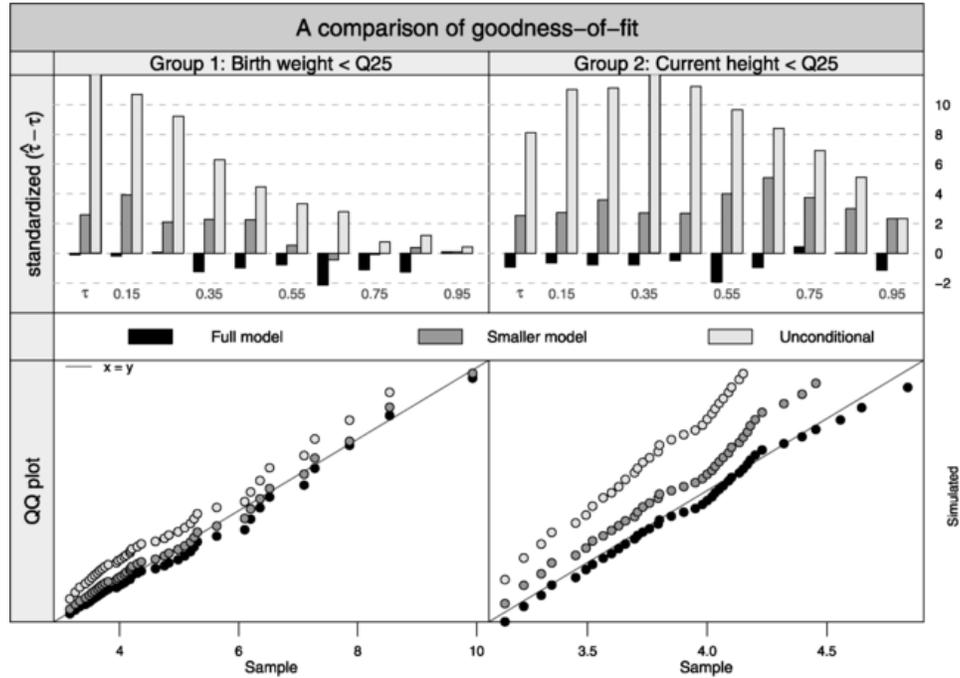

FIG. 4. *Model assessment for Groups* 1 *and* 2.

**7. Further discussion.** It is clear that the conditional growth charts that we propose in this paper are not a set of charts in the usual sense. Screening based on the global model on longitudinal growth data needs to be done with customized charts for each subject. However, such customized charts can be generated easily given the estimated models, and all one needs is to input the information on a given subject under screening, making the clinical use of the conditional growth charts really practical on any desktop.

In this section we discuss some limitations of and generalizations to the model (2.1). We hope that our work will stimulate future research on longitudinal growth charts.

7.1. *Generalization.* Fixed measurement intervals are usually assumed in an autoregressive model. By assuming that the AR coefficients are linear in the time spacings, the global model (2.1) has some flexibility to deal with irregular measurement times. We can further generalize the model as

$$(7.1) \qquad Y_{i,j} = g_0(t_{i,j}) + \sum_{k=1}^{p} b_k(D_{i,j,k}) Y_{i,j-k} + X_{i,j}^\top \gamma_0 + e_{i,j},$$

where the AR coefficients are some unknown functions of measurement time spacings. In fact, for small measurement time spacings, $a_k + b_k D_{i,j,k}$ works reasonably well as a linear approximation to the general function $b_k(D_{i,j,k})$. The estimate of $b_k(D_{i,j,k})$ in model (7.1) can be obtained by either a parametric or a nonparametric approach. For a parametric approach, we may specify a parametric coefficient function, such as $a_k + b_k D_{i,j,k} + c_k D_{i,j,k}^2$. To choose between $a_k + b_k D_{i,j,k} + c_k D_{i,j,k}^2$ and $a_k + b_k D_{i,j,k}$, we can simply test the hypothesis $H_0 : c_k = 0$. For a nonparametric approach, one may use the varying-coefficient models explored by Kim [25] or the dynamic models of Cai and Xu [2].

Another restrictive feature of the global model is that the AR coefficients do not depend on the measurement times, but only on their spacings. An alternative model can be set up as

$$(7.2) \quad Y_{i,j} = g_0(t_{i,j}) + \sum_{k=1}^{p} [a_k(t_{i,j-k}) + b_k(t_{i,j-k}) D_{i,j,k}] Y_{i,j-k} + X_{i,j}^\top \gamma_0 + e_{i,j},$$

where we allow the parameters $a_k$ and $b_k$ to be functions of prior measurement times $t_{i,j-k}$. Time series models studied by Xu [44] may also be considered. In practice, the flexibility of the models becomes an advantage only when the sample size is large.

7.2. *Localization.* If all the subjects have measurements on a series of fixed time points, say $(t_0, t_1, t_2, \ldots, t_p)$, then $a_k$ and $b_k$ in the global model are no longer identifiable. To model the conditional quantiles of the measurement at time $t_0$ given the prior path and other covariates, one may use

$$(7.3) \qquad Y_{t_0,i} = \beta_0 + \sum_{k=1}^{p} \beta_k Y_{t_k,i} + \gamma^\top X_i + e_i, \qquad i = 1,\ldots,n,$$

where $Y_{t,i}$ is the measurement of the $i$th subject at time $t$, $(t_1, t_2, \ldots, t_p)$ are $p$ prespecified earlier time points, which are not necessarily evenly spaced, $X_i$ is an $l$-dimensional covariate and $e_i$ are i.i.d. random variables with zero $\tau$th quantile. It follows that the conditional $\tau$th quantile of $Y_{t_0,i}$ is a linear combination of $p$ prior measurements $Y_{t_1,i}, Y_{t_2,i}, \ldots, Y_{t_p,i}$, and the other covariate $X_i$.

7.3. *Conditional versus marginal growth charts.* We have focused on the conditional growth charts in the present paper, because conditional charts are highly valuable in evaluating recent growth status of subjects. However, we do not suggest that we should do away with marginal growth charts, which are most commonly used by health professionals today.

The conditional growth charts may be unsuccessful in screening out subjects with gradual but persistent slowdown in growth. Further studies are needed to address issues related to the "longitudinal" use of growth charts over time. If a subject consistently falls below the 25th conditional percentile for several time periods, it might signal a growth problem. At the present time, the marginal growth charts may be used in conjunction with conditional charts.

7.4. *Monotonicity.* The growth charts from model (2.1) require estimation of a series of conditional quantile functions. It is not guaranteed that the estimated quantiles are monotone in $\tau$. For a discussion of constrained quantile regression that avoids this crossing problem, see [16].

## 8. Technical proofs.

8.1. *Proof of Theorem* 3.1. First we write model (2.1) in matrix form as

$$(8.1) \qquad Y_n = g_0(T) + \tilde{X}_n^\top \boldsymbol{\beta}_0 + e_n,$$

where $T = (t_{i,j})_{N \times 1} = (t_{1,2}, \ldots, t_{1,m_1}, t_{2,2}, \ldots, t_{2,m_2}, \ldots, t_{n,2}, \ldots, t_{n,m_n})^\top$, $Y_n = (Y_{i,j})_{N \times 1}$, $\tilde{X}_n = (\tilde{X}_{i,j}^\top)_{N \times (l+2)}$, $e_n = (e_{i,j})_{N \times 1}$ and $\boldsymbol{\beta}_0 = (a_0, b_0, \gamma_0^\top)^\top$. The indices $(i,j)$ in the vectors $Y_n$ and $e_n$ and the matrix $\tilde{X}_n$ are arranged in the same way as in $T$. Here $N = \sum_{i=1}^{n}(m_i - 1)$ is the total number of observations used for the estimation. Note that $m_i$ are uniformly bounded for all $n$, $O(N) = O(n)$. Recall that $Z_n = (\pi_{1,2}, \pi_{1,3}, \ldots, \pi_{i,j}, \ldots, \pi_{n,m_n})^\top_{N \times k_n}$ and $G_n = Z_n(Z_n^\top Z_n)Z_n^\top$. We further define

$$H_n^2 = k_n Z_n^\top Z_n, \qquad z_n = (z_{n1}, z_{n2}) = (k_n^{1/2} Z_n H_n^{-1}, n^{-1/2}(I - G_n)\tilde{X}_n),$$

$$\theta(\boldsymbol{\alpha}, \boldsymbol{\beta}) = (\theta_1(\boldsymbol{\alpha},\boldsymbol{\beta})^\top, \theta_2(\boldsymbol{\alpha},\boldsymbol{\beta})^\top)^\top = \begin{pmatrix} k_n^{-1/2} H_n \boldsymbol{\alpha} + k_n^{1/2} H_n^{-1} Z_n^\top \tilde{X}_n \boldsymbol{\beta} \\ n^{1/2} \boldsymbol{\beta} \end{pmatrix}.$$

Note that $z_n \neq Z_n$ in our notation; $Z_n$ is the matrix of the spline basis functions, and $z_n$ is the whole design matrix after normalization and orthogonalization. We split $z_n$ into two parts: the matrix $z_{n1}$ is the normalized design matrix of the spline bases, and $z_{n2}$ is the normalized design matrix of the linear components that is orthogonal to $z_{n1}$. Also, we denote $z_{n1,ij}$ and $z_{n2,ij}$ as the column components associated with the $i$th subject and the $j$th measurement. Similar partitions and standardizations were used in [21].

With this notation, model (8.1) can be further written as

$$\begin{aligned} (8.2) \quad Y_n &= z_{n1}\theta_1(\boldsymbol{\alpha}_0, \boldsymbol{\beta}_0) + z_{n2}\theta_2(\boldsymbol{\alpha}_0, \boldsymbol{\beta}_0) + R_n + e_n \\ &= z_n \theta(\boldsymbol{\alpha}_0, \boldsymbol{\beta}_0) + R_n + e_n. \end{aligned}$$

The orthogonality between $z_{n1}$ and $z_{n2}$ will not affect the estimation of $\boldsymbol{\alpha}_0$ and $\boldsymbol{\beta}_0$, but could simplify our proof. Using this orthogonalized model format, let

$$C_{i,j}(\theta) = \rho_\tau(e_{i,j} - z_{n,ij}\theta - R_{nij}) - \rho_\tau(e_{i,j} - R_{nij}) + z_{n,ij}\theta\psi(e_{i,j}),$$

$$\begin{aligned} D_{i,j}(\theta_1, \theta_2) &= \rho_\tau(e_{i,j} - z_{n1,ij}\theta_1 - z_{n2,ij}\theta_2 - R_{nij}) \\ &\quad - \rho_\tau(e_{i,j} - z_{n1,ij}\theta_1 + R_{nij}) + z_{n2,ij}\theta_2\psi(e_{i,j}). \end{aligned}$$

The proof of Theorem 3.1 relies on the following three lemmas.

LEMMA 8.1. *Under Assumptions 1–7, for any sequence $\{L_n\}$ such that $1 < L_n < k_n^{\delta_0}$ for some $\delta_0 \in (0, r/10)$, we have*

$$(8.3) \quad \sup_{\|\theta\| < L_n k_n^{1/2}} k_n^{-1} \left| \sum_{i=1}^n \sum_{j=2}^{m_i} (C_{i,j}(\theta) - E[C_{i,j}(\theta)]) \right| = o_p(1),$$

*and for any $M_1 > 0$ and $M_2 > 0$, we have*

$$(8.4) \quad \sup_{\|\theta_1\| \leq M_1 k_n^{1/2}; \|\theta_2\| \leq M_2} \left| \sum_{i=1}^n \sum_{j=2}^{m_i} (D_{i,j}(\theta_1, \theta_2) - E[D_{i,j}(\theta_1, \theta_2)]) \right| = o_p(1).$$

PROOF. The uniform convergence of (8.3) and (8.4) can be proven following arguments similar to those used in [19]. We skip the details to save space. □

LEMMA 8.2. *Under the same conditions of Lemma* 8.1, *we have, as* $n \to \infty$,

$$P\left(\inf_{|\theta|=1} k_n^{-1} \left|\sum_{i=1}^{n}\sum_{j=2}^{m_i}(E[C_{i,j}(L_n k_n^{1/2}\theta)] \right.\right.$$
(8.5)
$$\left.\left. - L_n k_n^{1/2} z_{n,ij}\theta\psi(e_{i,j}))\right| > 1\right) \to 1.$$

PROOF. First, under the condition of $\|\theta\| = 1$ and $d_n < W_5$, it follows from Assumption 5 that $E[C_{i,j}(\theta)] = E[\int_{-R_{nij}}^{-L_n k_n^{1/2} z_{n,ij}\theta - R_{nij}} (b_{i,j}s + B_n s^2)\,ds]$, where $B_n$ is a bounded sequence. Let $\tilde{d}_n = \max_{i,j}(L_n k_n^{1/2}\|z_{n,ij}\| + |R_{nij}|)$; we have

$$k_n^{-1}\sum_{i=1}^{n}\sum_{j=2}^{m_i} E[C_{i,j}(L_n k_n^{1/2}\theta)]$$

$$\geq k_n^{-1}\sum_{i=1}^{n}\sum_{j=2}^{m_i} E\left[\frac{b_{i,j}}{2}(-L_n k_n^{1/2} z_{n,ij}\theta - R_{nij})^2 - \frac{b_{i,j}}{2}(R_{nij})^2\right]$$

(8.6)
$$- E\left[B_n \tilde{d}_n^2 k_n^{-1}\sum_{i=1}^{n}\sum_{j=2}^{m_i} |L_n k_n^{1/2} z_{n,ij}\theta|\right]$$

$$= \left[L_n^2/2\sum_{i=1}^{n}\sum_{j=2}^{m_i} E[b_{i,j}(z_{n,ij}\theta)^2]\right.$$

$$\left. + L_n k_n^{-1/2}\sum_{i=1}^{n}\sum_{j=2}^{m_i} E(b_{i,j}R_{nij}z_{n,ij}\theta)\right] + o(1).$$

To see the last step of (8.6), we first note that $E\|z_{n,ij}\| = O(n^{-1/2}k_n^{1/2})$, $\|\theta\| = 1$, and thus $k_n^{1/2}\sum_{i=1}^{n}\sum_{j=2}^{m_i} E|z_{n,ij}\theta| < \infty$. Meanwhile, $\tilde{d}_n^2 L_n \leq d_n = o_p(1)$. It follows from the dominated convergence theorem that $E[\tilde{d}_n^2 L_n k_n^{-1/2} \times \sum_{i=1}^{n}\sum_{j=2}^{m_i} |z_{n,ij}\theta|] = o(1)$, and therefore the last step of (8.6) holds.

Since $b_{i,j}$ is uniformly bounded away from zero and infinity due to Assumption 5, it is easy to see that $\sum_{i=1}^{n}\sum_{j=2}^{m_i} E[b_{i,j}(z_{n,ij}\theta)^2]$ is bounded away from 0 and from infinity, uniformly in $\|\theta\| = 1$. On the other hand, we have $k_n^{-1/2}\sum_{i=1}^{n}\sum_{j=2}^{m_i} E(b_{i,j}R_{nij}z_{n,ij}\theta) = O(L_n)$. Therefore, we conclude that

(8.7)
$$k_n^{-1}\sum_{i=1}^{n}\sum_{j=2}^{m_i} E[C_{i,j}(L_n k_n^{1/2}\theta)] \geq cL_n^2,$$

for some $c > 0$, which implies that we only need to show

$$(8.8) \quad L_n k_n^{-1/2} \sum_{i=1}^{n} \sum_{j=2}^{m_i} z_{n,ij} \theta \psi(e_{i,j}) = o_p(L_n^2)$$

to prove Lemma 8.2.

By Theorem 4 of [3], the eigenvalues of $k_n Z_n^\top Z_n / n$ are bounded away from zero and infinity when $n$ is sufficiently large. Assumption 4 implies that the eigenvalues of $n^{-1} \psi_\tau(e_n) \psi_\tau(e_n)^\top$ are also bounded. Based on these facts, if we denote $\psi_\tau(e_n) = (\psi_\tau(e_{1,1}), \ldots, \psi_\tau(e_{n,m_n}))^\top$, then $\psi_\tau(e_n)^\top \psi_\tau(e_n) = \mathrm{diag}(\psi_\tau(e_i)^\top \psi_\tau(e_i))$. Let $\lambda_i$ be the largest eigenvalue of $\psi_\tau(e_i)^\top \psi_\tau(e_i)$; $\sup_i \lambda_i = O_p(1)$ by Assumption 4. It thus follows that

$$(8.9) \quad \begin{aligned} \sup_{\|\theta\|=1} L_n k_n^{-1/2} \sum_{i=1}^{n} \sum_{j=2}^{m_i} z_{n,ij} \theta \psi_\tau(e_{i,j}) \\ = \sup_{\|\theta\|=1} L_n k_n^{-1/2} [\theta^\top z_n^\top \psi_\tau(e_n) \psi_\tau(e_n)^\top z_n \theta]^{1/2} \\ \leq \sup_{\|\theta\|=1} \sup_i \lambda_i L_n k_n^{-1/2} [\theta^\top z_n^\top z_n \theta]^{1/2} \\ = \left(\sup_i \lambda_i\right) L_n k_n^{-1/2} = O_p(L_n k_n^{-1/2}). \end{aligned}$$

Therefore, (8.8) holds, and the proof of Lemma 8.2 is complete. □

LEMMA 8.3. *Under Assumptions 1–9, we have*

$$(8.10) \quad \sup_{\|\theta_1\| \leq L k_n^{1/2}; \|\theta_2\| \leq M} \left| \sum_{i=1}^{n} \sum_{j=2}^{m_i} E[D_{i,j}(\theta_1, \theta_2)] - n^{-1} \tfrac{1}{2} \theta_2^\top K_n \theta_2 \right| = o_p(1).$$

PROOF. First of all, we note that

$$\sup_{\substack{\|\theta_1\| \leq L k_n^{1/2}; \\ \|\theta_2\| \leq M}} \left| \sum_{i=1}^{n} \sum_{j=2}^{m_i} E[D_{i,j}(\theta_1, \theta_2)] - \tfrac{1}{2} n^{-1} \theta_2^\top K_n \theta_2 \right|$$

$$= \sup_{\substack{\|\theta_1\| \leq L k_n^{1/2}; \\ \|\theta_2\| \leq M}} \left| \sum_{i=1}^{n} \sum_{j=2}^{m_i} E \left\{ E_e \left[ \int_{-z_{n1,ij}\theta_1 - R_{nij}}^{-z_{n1,ij}\theta_1 - z_{n2,ij}\theta_2 - R_{nij}} \psi(e_{i,j} + s) \, ds \right. \right.\right.$$

$$(8.11) \quad \left.\left.\left. |z_{n,ij}| \right] \right\} - n^{-1} \tfrac{1}{2} \theta_2^\top K_n \theta_2 \right|$$

$$= \sup_{|\theta_2\| \leq M} |\tfrac{1}{2}[E(\theta_2^\top z_{n2}^\top B_n z_{n2} \theta_2) - \theta_2^\top z_{n2}^\top B_n z_{n2} \theta_2] + E(\theta_2^\top z_{n2} r_n)|$$

$$+ o(1),$$

where $r_n$ is the $n \times 1$ vector with elements $(b_{i,j} R_{nij})$.

Assumption 8 implies that

$$\sup_{\|\theta_2\| \leq M} |E(\theta_2^\top z_{n2}^\top B_n z_{n2} \theta_2) - \theta_2^\top z_{n2}^\top B_n z_{n2} \theta_2| = o_p(1). \tag{8.12}$$

On the other hand, $\sup_{\|\theta_2\| \leq M} \theta_2^\top z_{n2}^\top r_n \leq k_n^{-r} \sup_{\|\theta_2\| \leq M} \theta_2^\top z_{n2}^\top \mathbf{1}_n$, and it follows that

$$\sup_{\|\theta_2\| \leq M} E[\theta_2^\top z_{n2}^\top r_n] = o(1). \tag{8.13}$$

Combining (8.11)–(8.13), we arrive at (8.10). Lemma 8.3 is thus proven. □

With Lemmas 8.1 to 8.3, the proof of Theorem 3.1 can be outlined as follows:

PROOF OF THEOREM 3.1. Model (8.1) can be reconstructed as

$$Y_n = z_{n1} \theta_1(\boldsymbol{\alpha}_0, \boldsymbol{\beta}_0) + z_{n2} \theta_2(\boldsymbol{\alpha}_0, \boldsymbol{\beta}_0) + R_n + e_n$$
$$= z_n \theta(\boldsymbol{\alpha}_0, \boldsymbol{\beta}_0) + R_n + e_n.$$

We define

$$\tilde{\theta}_n = (\tilde{\theta}_{n1}^\top, \tilde{\theta}_{n2}^\top)^\top$$
$$= (\theta_1(\hat{\boldsymbol{\alpha}}_n, \hat{\boldsymbol{\beta}}_n) - \theta_1(\boldsymbol{\alpha}_0, \boldsymbol{\beta}_0), \theta_2(\hat{\boldsymbol{\alpha}}_n, \hat{\boldsymbol{\beta}}_n) - \theta_2(\boldsymbol{\alpha}_0, \boldsymbol{\beta}_0))^\top. \tag{8.14}$$

The objective function (2.2) can be written as

$$\sum_{i=1}^{n} \sum_{j=2}^{m_i} \rho_\tau(e_{i,j} - z_{n,ij} \theta_n - R_{nij}), \tag{8.15}$$

and we have the fact that $\hat{\boldsymbol{\beta}}_n$ minimizes (2.2) and $\tilde{\theta}_n$ minimizes (8.15).

Directly from the definition of $\tilde{\theta}_n$, we have

$$\|\hat{\boldsymbol{\beta}}_n - \boldsymbol{\beta}_0\| = \|n^{-1/2} \tilde{\theta}_{n2}\| = O_p(n^{-1/2} \|\tilde{\theta}_{n2}\|), \tag{8.16}$$

which implies that $\|\tilde{\theta}_{n2}\| = O_p(k_n^{-1/2})$ is a sufficient condition for (3.1). On the other hand,

$$\frac{1}{n} \sum_{i=1}^{n} \sum_{j=2}^{m_i} (\hat{g}_n(t_{i,j}) - g_0(t_{i,j}))^2$$
$$\leq \frac{2}{n} \sum_{i=1}^{n} \sum_{j=2}^{m_i} (\pi_{i,j}(\hat{\boldsymbol{\alpha}}_n - \boldsymbol{\alpha}_0))^2 + 2W_3^2 k_n^{-2r}$$
$$\leq 2[n^{-1} \|\tilde{\theta}_{n1}\|^2 + n^{-1} \|H_n^{-1} k_n^{1/2} Z_n^\top \tilde{X}_n(\hat{\boldsymbol{\beta}}_n - \boldsymbol{\beta}_0)\|^2] + 2W_3^2 k_n^{-2r}$$
$$= O_p(n^{-1} \|\tilde{\theta}_{n1}\|^2) + O(\|\hat{\beta} - \beta_0\|) + O(k_n^{-2r}). \tag{8.17}$$

The last step of (8.17) is due to the fact that $n^{-1}\|G_n\tilde{X}_n\|^2 \leq n^{-1}\|\tilde{X}_n\|^2 < \infty$ under Assumption 7. Combining (8.16) and (8.17), we conclude that it is sufficient to show $\|\tilde{\theta}_n\| = O_p(k_n^{-1/2})$.

According to Lemmas 8.1 and 8.2, for any $\varepsilon$ there exists $L_\varepsilon$ such that

$$P\left\{\inf_{\|\theta\|>L_\varepsilon k_n^{1/2}} \sum_{i=1}^n \sum_{j=2}^{m_i} \rho_\tau(e_{i,j} - z_{n,ij}\theta - R_{nij}) > \sum_{i=1}^n \sum_{j=2}^{m_i} \rho_\tau(e_{i,j} - R_{nij})\right\} > 1 - \varepsilon.$$

Since $\tilde{\theta}_n$ minimizes $\sum_{i=1}^n \sum_{j=2}^{m_i} \rho_\tau(e_{i,j} - z_{n,ij}\theta - R_{nij})$ over the space $\mathbf{R}^{p_n}$, we have $P(\|\tilde{\theta}_n\| < L_\varepsilon k_n^{1/2}) > 1 - \varepsilon$, and thus $\|\tilde{\theta}_n\| = O_p(k_n^{1/2})$.

With the consistency of the global model, we can take a further step to show the asymptotic normality of the coefficient estimate $\hat{\beta}_n$.

We denote by $z_{n2}^{(i)}$ the submatrix of $z_{n2}$ corresponding to the $i$th subject; in other words, $z_{n2} = (z_{n2}^{(1)T}, z_{n2}^{(2)T}, \ldots, z_{n2}^{(n)T})^\top$. Let $\theta_{n2}^* = nK_n^{-1}\sum_{i=1}^n z_{n2}^{(i)T}\psi(e_i)$. Due to Assumptions 8 and 9, $\theta_{n2}^*$ is asymptotically normally distributed with asymptotic variance–covariance matrix $K^{-1}SK^{-1}$. Since $\tilde{\theta}_{n2} = n^{1/2}\hat{\beta}_n$, if we can show that

(8.18) $$\|\tilde{\theta}_{n2} - \theta_{n2}^*\| = o_p(1),$$

then (3.3) holds.

Due to the definition of $\theta_{n2}^*$ and the consistency of $\tilde{\theta}_n$, we know that $P(\theta_{n2}^* < M) \to 1$ and $P(\|\tilde{\theta}_{n1}\| < Lk_n^{1/2}) \to 1$ for any $L > 0$ and $M > 0$. Let

$$\tilde{D}_{i,j}(\theta_2, \theta_2^*) = \rho_\tau(e_{i,j} - z_{n1,ij}\tilde{\theta}_{n1} - z_{n2,ij}\theta_2 - R_{nij})$$
$$- \rho_\tau(e_{i,j} - z_{n1,ij}\tilde{\theta}_{n1} - z_{n2,ij}\theta_2^* - R_{nij}).$$

It follows from Lemma 8.4 that, for any given $\delta > 0$,

$$\sup_{\|\theta_2-\theta_{2n}^*\|\leq\delta}\left|\sum_{i=1}^n\sum_{j=2}^{m_i}\{\tilde{D}_{i,j}(\theta_2,\theta_{n2}^*) - z_{n2,ij}\psi(e_{i,j})(\theta_2-\theta_{n2}^*)\right.$$
$$\left. - E[\tilde{D}_{i,j}(\theta_2,\theta_{n2}^*)]\}\right| = o_p(1).$$

Furthermore, by Lemma 8.3 we have

$$\sup_{|\theta_2-\theta_{2n}^*|\leq\delta}\left|\sum_{i=1}^n\sum_{j=2}^{m_i}[\tilde{D}_{i,j}(\theta_2,\theta_{n2}^*)]\right.$$
$$\left.+ (\theta_2-\theta_{n2}^*)^\top\sum_{i=1}^n z_{n2}^{(i)T}\psi(e_i) - n^{-1}\tfrac{1}{2}\theta_2^\top K_n\theta_2 + n^{-1}\tfrac{1}{2}\theta_{2n}^{*\top}K_n\theta_{n2}^*\right|$$

(8.19) $$= \sup_{|\theta_2-\theta_{2n}^*|\leq\delta}\left|\sum_{i=1}^n\sum_{j=2}^{m_i}[\tilde{D}_{i,j}(\theta_2,\theta_{n2}^*)]\right.$$

$$+ n^{-1}[(\theta_2 - \theta_{n2}^*)^\top K_n \theta_2^* - \tfrac{1}{2}\theta_2^\top K_n \theta_2 + \tfrac{1}{2}\theta_{2n}^{*\top} K_n \theta_{n2}^*] \Bigg|$$

$$= \sup_{|\theta_2 - \theta_{2n}^*| \leq \delta} \left| \sum_{i=1}^n \sum_{j=2}^{m_i} [\tilde{D}_{i,j}(\theta_2, \theta_{n2}^*)] - n^{-1}\tfrac{1}{2}(\theta_2 - \theta_{n2}^*)^\top K_n(\theta_2 - \theta_{n2}^*) \right|$$

$$= o_p(1).$$

For sufficiently large $n$, $n^{-1}(\theta_2 - \theta_{n2}^*)^\top K_n(\theta_2 - \theta_{n2}^*) > 0$ when $\|\theta_2 - \theta_{n2}^*\| > \delta$. Then (8.19) implies that

$$\lim_{n \to \infty} P\Bigg( \inf_{\|\theta_2 - \theta_2^*\| \geq \delta} \sum_{i=1}^n \sum_{j=2}^{m_i} \rho_\tau(e_{i,j} - z_{n1,ij}\tilde{\theta}_{n1} - z_{n1,ij}\theta_2 - R_{nij})$$

(8.20)
$$> \sum_{i=1}^n \sum_{j=2}^{m_i} \rho_\tau(e_{i,j} - z_{n1,ij}\tilde{\theta}_{n1} - z_{n1,ij}\theta_2^* - R_{nij}) \Bigg) = 1.$$

Since $\tilde{\theta}_{n1}$ minimizes $\sum_{i=1}^n \sum_{j=2}^{m_i} \rho_\tau(e_{i,j} - z_{n1,ij}\theta_1 - z_{n1,ij}\theta_2 - R_{nij})$ over $\mathbf{R}^{l+2}$, (8.20) implies that for any $\delta$, $P(\|\tilde{\theta}_{n2} - \theta_{n2}^*\| > \delta) \to 0$, that is, $\|\tilde{\theta}_{n2} - \theta_{n2}^*\| = o_p(1)$. The proof of Theorem 3.1 is hence complete. $\square$

8.2. *Proof of Theorem* 4.1. Recall that $S_n = n^{-1/2}\sum_{i=1}^n \sum_{j=2}^{m_i} \psi_\tau(Y_{i,j} - w_{i,j}^\top \hat{\phi}_n)v_{i,j}$. Since the $v_{i,j}$'s are the least squares residuals from regressing $\tilde{X}_{n1}$ on $\mathcal{W}_n$, they differ from $v_{i,j}^{(0)} = \tilde{X}_{n1,ij} - E\{\tilde{X}_{n1,ij}|\mathcal{W}_{n,ij}\}$ by $O_p(k_n^{-r})$. It is easy to show that the limiting distribution of $T_n = S_n^\top V_n^{-1} S_n$ will not change if $v_{i,j}$ are replaced by $v_{i,j}^{(0)}$; the latter enjoy intersubject independence and are often easier to handle mathematically. To simplify notation, we will simply prove the results in this section by assuming that the $v_{i,j}$'s are intersubject independent.

First, we give the following two lemmas.

LEMMA 8.4. *Let* $u_{i,j}(\phi, \phi_0) = \varphi_\tau(Y_{i,j} - w_{i,j}^\top \phi)v_{i,j} - \varphi_\tau(Y_{i,j} - w_{i,j}^\top \phi_0)v_{i,j} - E\varphi_\tau(Y_{i,j} - w_{i,j}^\top \phi)v_{i,j} + E\varphi_\tau(Y_{i,j} - w_{i,j}^\top \phi_0)v_{i,j}$. *Under the assumptions of Theorem* 4.1, *for any* $L > 0$, *we have*

$$\sup_{\|\phi - \phi_0\| \leq L(k_n/n)^{1/2}} \left\| \sum_{i=1}^n \sum_{j=2}^{m_i} u_{i,j}(\phi, \phi_0) \right\| = O_p(n^{1/4} k_n^{3/4}(\ln(n))^{1/2})$$

(8.21)
$$= o_p(n^{-1/2}).$$

PROOF. Lemma 8.4 can be viewed as a special case of He and Shao [17], which considered the asymptotic behavior of $M$-estimators with increasing dimension.

Since $\varphi_\tau(t) = 1 - \tau I_{\{t<0\}}$ is a constant function with a jump at $t=0$, we have

$$[\varphi_\tau(Y_{i,j} - w_{i,j}^\top \phi) - \varphi_\tau(Y_{i,j} - w_{i,j}^\top \phi_0)]^2 \|v_{i,j}\|^2 \tag{8.22}$$
$$\leq 4\|v_{i,j}\|^2 I_{\{|Y_{i,j} - w_{i,j}^\top \phi| \leq |w_{i,j}^\top(\phi-\phi_0)|\}}.$$

Under Assumption D3,

$$E \sum_{i=1}^n \left\| \sum_{j=2}^{m_i} u_{i,j}(\phi, \phi_0) \right\|^2 \leq 8\kappa \left(\sup_i m_i\right) E \sum_{i=1}^n \sum_{j=2}^{m_i} \|v_{i,j}\|^2 |w_{i,j}^\top(\phi-\phi_0)|, \tag{8.23}$$

which, together with Assumption D2, implies

$$A_n = \sup_{\|\phi-\phi_0\| \leq L(k_n/n)^{1/2}} \sum_{i=1}^n E_\phi \left\| \sum_{j=2}^{m_i} u_{i,j}(\phi, \phi_0) \right\|^2 = O((nk_n)^{1/2}). \tag{8.24}$$

Moreover, since $\varphi_\tau(t)$ is bounded and $\max_{i,j} \|v_{i,j}\|^2 = O_p(1)$, it follows from Lemma 2.2 of [17] that

$$B_n = \sup_{\|\phi-\phi_0\| \leq L(k_n/n)^{1/2}} \sum_{i=1}^n \left\| \sum_{j=2}^{m_i} u_{i,j}(\phi, \phi_0) \right\|^2 = O_p((nk_n)^{1/2}). \tag{8.25}$$

Finally, combining (8.23) and Assumption D2, we note that condition (C1) of Lemma 3.3 of [17] is satisfied, and (8.21) holds consequently, where the right-hand side comes from $O_p((\ln(n)k_n)^{1/2}(n^{-2} + A_n^{1/2} + B_n^{1/2}))$. □

LEMMA 8.5. *Under the assumptions of Theorem* 4.1,

$$\sup_{\|\phi-\phi_0\| \leq L(k_n/n)^{1/2}} n^{-1/2} \left\| \sum_{i=1}^n \sum_{j=2}^{m_i} E\varphi_\tau(Y_{i,j} - w_{i,j}^\top \phi) v_{i,j} \right.$$
$$\left. - E\varphi_\tau(Y_{i,j} - w_{i,j}^\top \phi_0) v_{i,j} \right\| = o(1). \tag{8.26}$$

PROOF. Under the condition $\|\phi - \phi_0\| < L(k_n/n)^{1/2}$, we first expand the conditional mean $E[\varphi_\tau(Y_{i,j} - w_{i,j}^\top \phi) v_{i,j} | v_{i,j}, w_{i,j}]$ around $\phi_0$ for each $i$ and $j$. Note that, for a real-valued random variable $u$ $E[\varphi_\tau(u)] = \tau - F_u(0)$, where $F_u$ is the distribution of $u$, we have

$$E\{[\varphi_\tau(Y_{i,j} - w_{i,j}^\top \phi) - \varphi_\tau(Y_{i,j} - w_{i,j}^\top \phi_0)] v_{i,j} | v_{i,j}, w_{i,j}\}$$
$$= f_{e_{i,j}}(0) v_{i,j} w_{i,j}^\top (\phi-\phi_0) \tag{8.27}$$
$$+ \tfrac{1}{2} f'_{e_{i,j}}(0) v_{i,j} (\phi-\phi_0)^\top w_{i,j} w_{i,j}^\top (\phi-\phi_0) + o(n^{-1}).$$

The remainder term in (8.27) is $o(n^{-1})$ because $\|\phi - \phi_0\| < L(k_n/n)^{1/2}$ and $k_n = o(n^{1/4})$. It follows from (8.27) that the left-hand side of (8.26) is

$$\sup_{\|\phi-\phi_0\|\leq L(k_n/n)^{1/2}} n^{-1/2} \left\| \sum_{i=1}^{n}\sum_{j=2}^{m_i} E\varphi_\tau(Y_{i,j} - w_{i,j}^\top \phi)v_{i,j} \right.$$
$$\left. - E\varphi_\tau(Y_{i,j} - w_{i,j}^\top \phi_0)v_{i,j} \right\|$$

$$(8.28) \quad = \sup_{\|\phi-\phi_0\|\leq L(k_n/n)^{1/2}} n^{-1/2} \| bE\mathcal{V}_n^\top \mathcal{W}_n(\phi - \phi_0)$$
$$+ cE(\phi - \phi_0)^\top \mathcal{W}_n[\text{diag}(v_{i,j})]\mathcal{W}_n^\top (\phi - \phi_0) \|$$
$$+ o(1),$$

where $b = f_{e_{i,j}}(0)$ and $c = f'_{e_{i,j}}(0)$. Note that $\mathcal{V}_n$ and $\mathcal{W}_n$ are orthogonal to each other, that is, $\mathcal{V}_n^\top \mathcal{W}_n = 0$, the first term in (8.28) is zero. Moreover,

$$\sup_{\|\phi-\phi_0\|\leq L(k_n/n)^{1/2}} |E(\phi - \phi_0)^\top \mathcal{W}_n[\text{diag}(v_{i,j})]\mathcal{W}_n^\top (\phi - \phi_0)|$$
$$\leq L^2(k_n/n) E\left[\max_{i,j} v_{i,j} \|\mathcal{W}_n\|^2\right] = O(k_n^2).$$

The last equation holds due to Assumption D2 and the fact that $E\|\mathcal{W}_n\|^2 = O(nk_n)$. When $k_n = o(n^{1/4})$,

$$(8.29) \quad n^{-1/2} cE(\phi - \phi_0)^\top \mathcal{W}_n \text{diag}(v_{i,j})\mathcal{W}_n^\top (\phi - \phi_0) = o(1).$$

Thus, (8.26) follows from (8.28) and (8.29). □

PROOF OF THEOREM 4.1. In Theorem 3.1 we have shown the consistency of $\hat{\phi}_n$ given $k_n \approx n^{1/(2r+1)}$, which, however, is not the necessary condition for the consistency. In fact, for any $k_n = o(n^{1/3})$, the consistency of $\hat{\phi}_n$ holds at the convergence rate $(k_n/n)^{1/2}$, that is, $\|\hat{\phi}_n - \phi_0\| = O_p(n^{-1/2}k_n^{1/2})$. Therefore, we can derive from Lemmas 8.4 and 8.5 that, when $k_n = o(n^{1/4})$,

$$(8.30) \quad n^{-1/2} \left\| \sum_{i=1}^{n}\sum_{j=2}^{m_i} \{\varphi_\tau(Y_{i,j} - w_{i,j}^\top \hat{\phi}_n)v_{i,j} \right.$$
$$\left. - \varphi_\tau(Y_{i,j} - w_{i,j}^\top \phi_0)v_{i,j}\} \right\| = o_p(1).$$

Denote $S_n^* = \sum_{i=1}^{n}\sum_{j=2}^{m_i}[\varphi_\tau(e_{i,j})v_{i,j}]$; then the summands $\sum_{j=2}^{m_i}[\varphi_\tau(e_{i,j})v_{i,j}]$ are independent of each other and have mean zero. Due to the between-

CONDITIONAL GROWTH CHARTS 27

subject independence,

$$\text{Var}(S_n^*) = \sum_{i=1}^{n} \text{Var}\left(\sum_{j=2}^{m_i}[\varphi_\tau(e_{i,j})v_{i,j}]\right) = \sum_{i=1}^{n} \mathcal{V}_i^\top \text{Var}(\varphi_\tau(e_i))\mathcal{V}_i$$

$$= \tau(1-\tau)\sum_{i=1}^{n} \mathcal{V}_i^\top \mathcal{V}_i = \tau(1-\tau)\tilde{X}_{n1}(I_n - G)\tilde{X}_{n1}^\top,$$

where $\mathcal{V}_i$ is the $m_i \times m_i$ submatrix of $\mathcal{V}_n$ corresponding to the $i$th subject, and $\varphi_\tau(e_i)$ is the same as defined at the beginning of Section 3. In the i.i.d. case, $\text{Var}(\varphi_\tau(e_i)) = \tau(1-\tau)I_{m_i}$. Let $V_n = n^{-1}\tau(1-\tau)\tilde{X}_{n1}(I_n - G)\tilde{X}_{n1}^\top$; it follows from the CLT that

$$(8.31) \qquad n^{-1/2}\sum_{i=1}^{n}\sum_{j=2}^{m_i}[\varphi_\tau(e_{i,j})v_{i,j}] \text{ is } AN(0, V_n).$$

Combining (8.30) and (8.31), it is clear that all we need to show Theorem 4.1 is

$$(8.32) \qquad n^{-1/2}\left\|\sum_{i=1}^{n}\sum_{j=2}^{m_i}[\varphi_\tau(Y_{i,j} - w_{i,j}^\top\hat{\phi}_n)v_{i,j} - \varphi_\tau(e_{i,j})v_{i,j}]\right\| = o_p(1).$$

However, it is hard to show (8.32) directly for some $k_n$, more specifically, for $n^{1/4r} \ll k_n \ll n^{1/2r}$ including $k_n \approx n^{1/2r+1}$. Instead, we take an intermediate step.

Let $K_n$ be the set of knots used in our global model estimation. We assume that the dimension of $K_n \approx n^{1/(2r+1)}$; the same proof goes through for other values of $K_n$. Under $H_0$, let $Y_{i,j} - w_{i,j}^\top\phi_0 = R_{nij} + e_{i,j}$, where $R_{nij}$ is the bias from the B-spline approximation using the set of knots $K_n$. The bias will go to zero with $\sup_{i,j}|R_{nij}| \leq W_3 k_n^{-r}$, where $W_3$ is defined in Assumption 1. By adding more knots into $K_n$, we have a new set of knots $\tilde{K}_n$ with its dimension, denoted as $\tilde{k}_n$, approximately $n^\alpha$, where $1/(2r) < \alpha < 1/4$. We also assume that the knots are added in such a way that the extended set of knots is quasiuniform. Using the new set of knots, we define $\tilde{\mathcal{W}}_n$, $\tilde{w}_{i,j}$, $\tilde{\mathcal{V}}_n$, $\tilde{v}_{i,j}$, $\tilde{\phi}$, $\tilde{\phi}_0$ and $\tilde{R}_{nij}$ the same way as $\mathcal{W}_n$, $w_{i,j}$, $\mathcal{V}_n$, $v_{i,j}$, $\hat{\phi}_n$, $\phi_0$ and $R_{nij}$ in the original setting.

Note that $v_{i,j}$ is the residual of a spline estimate for the $(i,j)$th element of $E[\tilde{X}_{n1}|\mathcal{W}_n]$, using $K_n$ as the set of knots; while $\tilde{v}_{i,j}$ is the residual of a spline estimate for the $(i,j)$th element of $E[\tilde{X}_{n1}|\tilde{\mathcal{W}}_n]$, using $\tilde{K}_n$ as the set of knots. Therefore, by the same arguments used for Theorem 3.1, $\frac{1}{n}\sum_i\sum_j \|v_{i,j} - \tilde{v}_{i,j}\|^2 = O_p(k_n^{-2r}) + O_p(\tilde{k}_n/n)$. In the meantime, by definition, $\tilde{\mathcal{V}}_n$ is orthogonal to both $\tilde{\mathcal{W}}_n$ and $\mathcal{W}_n$. Based on the facts above, following the similar arguments used for (8.30), we can show that, for $n^{1/2r} \ll \tilde{k}_n \ll n^{1/4}$ and $r > 2$,

$$(8.33) \quad n^{-1/2}\left\|\sum_{i=1}^{n}\sum_{j=2}^{m_i}\{\varphi_\tau(Y_{i,j} - w_{i,j}^\top\hat{\phi}_n)v_{i,j} - \varphi_\tau(Y_{i,j} - \tilde{w}_{i,j}^\top\tilde{\phi})\tilde{v}_{i,j}\}\right\| = o_p(1)$$

and

$$(8.34) \quad n^{-1/2}\left\|\sum_{i=1}^{n}\sum_{j=2}^{m_i}\{\varphi_\tau(Y_{i,j} - \tilde{w}_{i,j}^\top\tilde{\phi})\tilde{v}_{i,j} - \varphi_\tau(Y_{i,j} - \tilde{w}_{i,j}^\top\tilde{\phi}_0)\tilde{v}_{i,j}\}\right\| = o_p(1).$$

Moreover, noting that, given $v_{i,j}$, $E\varphi_\tau(Y_{i,j} - \tilde{w}_{i,j}^\top\tilde{\phi}_0) = E\varphi_\tau(\tilde{R}_{nij} + e_{i,j}) = O(\tilde{R}_{nij})$, we also have

$$(8.35) \quad n^{-1/2}\left\|\sum_{i=1}^{n}\sum_{j=2}^{m_i}\{\varphi_\tau(Y_{i,j} - \tilde{w}_{i,j}^\top\tilde{\phi}_0)v_{i,j} - \varphi_\tau(Y_{i,j} - \tilde{w}_{i,j}^\top\tilde{\phi}_0)\tilde{v}_{i,j}\}\right\| = o_p(1).$$

Combining (8.33)–(8.35), we have

$$(8.36) \quad n^{-1/2}\left\|\sum_{i=1}^{n}\sum_{j=2}^{m_i}\{\varphi_\tau(Y_{i,j} - w_{i,j}^\top\hat{\phi}_n)v_{i,j} - \varphi_\tau(Y_{i,j} - \tilde{w}_{i,j}^\top\tilde{\phi}_0)v_{i,j}\}\right\| = o_p(1).$$

Under Assumption D3, and the fact that $|\varphi_\tau(e_{i,j}+\tilde{R}_{nij})-\varphi_\tau(e_{i,j})| \leq I_{\{|e_{i,j}|\leq|\tilde{R}_{nij}|\}}$, we have

$$(8.37) \quad \begin{aligned} & n^{-1/2}E\sum_{i=1}^{n}\sum_{j=2}^{m_i}\|\varphi_\tau(Y_{i,j} - \tilde{w}_{i,j}^\top\tilde{\phi}_0)v_{i,j} - \varphi_\tau(e_{i,j})v_{i,j}\| \\ & \leq 2qn^{-1/2}E\sum_{i=1}^{n}\sum_{j=2}^{m_i}|\tilde{R}_{nij}|\|v_{i,j}\| = O(\tilde{k}_n^{-r}n^{1/2}) = o(1). \end{aligned}$$

The result (8.32) follows immediately from (8.30) and (8.37). The proof of Theorem 4.1 is thus complete. $\square$

**Acknowledgments.**The authors thank Anneli Pere for allowing us to use the Finnish growth data for illustration, and a Co-Editor and two referees for helpful comments and suggestions on our first draft.

DEPARTMENT OF BIOSTATISTICS
COLUMBIA UNIVERSITY
NEW YORK, NEW YORK 10032
USA
E-MAIL: yw2148@columbia.edu

DEPARTMENT OF STATISTICS
UNIVERSITY OF ILLINOIS
CHAMPAIGN, ILLINOIS 61820
USA
E-MAIL: x-he@uiuc.edu